# PING-PONG ON NEGATIVELY CURVED GROUPS

Rita Gitik

ABSTRACT. We generalize results about the classical Schottky groups to quasiconvex subgroups in negatively curved groups, answering a question of M. Bestvina.

## Introduction

The following result about classical Schottky groups of rank 2 goes back to Blaschke, Klein, Schottky and Poincaré. It can be viewed as a special case of the Klein-Maskit combination theorem ([Ma], p.135).

**Lemma.** *Let $h$ and $k$ be hyperbolic isometries of hyperbolic space such that their axes have disjoint endpoints. There exists a constant $C$ which depends only on $h$ and $k$, such that if $m > C$ and $n > C$, then the subgroup $<h^m, k^n>$ is isomorphic to the free product $<h^m> * <k^n>$.*

The well-known proof of this result can be abstracted into the following simple form. Let $H$ and $K$ be subgroups of a group $G$, and let $<H, K>$ be the smallest subgroup of $G$ containing $H$ and $K$.

**The Ping-Pong Lemma.** *Let $G$ be a group acting on a set $S$. Let $H$ and $K$ be subgroups of $G$ such that $|K : G_0| > 2$ or $|H : G_0| > 2$, where $G_0 = H \cap K$. If there exist disjoint non-empty subsets $S_K$ and $S_H$ of $S$ such that $K \setminus G_0$ maps $S_K$ into $S_H$ and $H \setminus G_0$ maps $S_H$ into $S_K$, then the subgroup $<H, K>$ is isomorphic to the amalgamated free product $<H, K> = H *_{G_0} K$.*

*Proof.* Let $h_1 k_1 \cdots h_n$ be an element of $<H, K> \setminus G_0$, where all $h_i$ and $k_i$ are not in $G_0$. Then for any $x \in S_H$, $(h_1 k_1 \cdots h_n)(x) \in S_K$, so $(h_1 k_1 \cdots h_n)(x) \neq x$, hence any element of odd syllable length in $<H, K>$ is non-trivial. Consider an element of even length $h_1 k_1 \cdots h_n k_n$. Without loss of generality, $|H : G_0| > 2$, hence there exists $h \in H$ such that $h \notin G_0$ and $h^{-1} h_1 \notin G_0$. But then the element $h^{-1} h_1 k_1 \cdots h_n k_n h$ has odd syllable length, hence it is non-trivial, therefore $h_1 k_1 \cdots h_n k_n$ is non-trivial.

Note that there exist counterexamples to the Lemma when $|K : G_0| = |H : G_0| = 2$, even if the action of $G$ on $S$ is free.

1991 *Mathematics Subject Classification.* 20F32, 20E06, 30F40.
Research partially supported by NSF grant DMS 9022140 at MSRI

Typeset by $\mathcal{A}_{\mathcal{M}}\mathcal{S}$-TEX

The Ping-Pong lemma allows us to describe the subgroup $<H, K>$. In general, even in the simple case when $H = <h>$ and $K = <k>$ are infinite cyclic groups with trivial intersection, we do not have much information about the group $<h, k>$. In particular, this group need not be isomorphic to the free group of rank two.

In this paper we describe the subgroup $<H, K>$ when $H$ and $K$ are quasiconvex subgroups of a negatively curved group $G$. We also give a condition for $<H, K>$ to be quasiconvex in $G$. In general, even if both $H = <h>$ and $K = <k>$ are infinite cyclic (hence quasiconvex in $G$, [Gr], p. 210), the group $<h, k>$ might be non-quasiconvex in $G$. This might happen even if $<h, k>$ is isomorphic to the free group of rank 2.

We prove the following results.

**Theorem 1.** *Let $H$ and $K$ be $\mu$-quasiconvex subgroups of a $\delta$-negatively curved group $G$. There exists a constant $C_0$, which depends only on $G, \delta$ and $\mu$, with the following property. For any subgroups $H_1 < H$ and $K_1 < K$ with $H_1 \cap K_1 = H \cap K = G_0$, if all the elements in $H_1$ and in $K_1$ which are shorter than $C_0$ belong to $G_0$, then $<H_1, K_1> = H_1 *_{G_0} K_1$. If, in addition, $H_1$ and $K_1$ are quasiconvex in $G$, then the subgroup $<H_1, K_1>$ is quasiconvex in $G$.*

Recall that a subgroup $H$ is malnormal in $G$ if for any $g \notin H$ the intersection of $H$ and $gHg^{-1}$ is trivial.

**Theorem 2.** *Let $H$ and $K$ be $\mu$-quasiconvex subgroups of a $\delta$-negatively curved group $G$. Assume that $H$ is malnormal in $G$. There exists a constant $C_1$, which depends only on $G, \delta$ and $\mu$, with the following property. For any subgroup $H_1 < H$ with $H_1 \cap K = H \cap K = G_0$, if all the elements in $H_1$ which are shorter than $C_1$ belong to $G_0$, then $<H_1, K> = H_1 *_{G_0} K$. If, in addition, $H_1$ is quasiconvex in $G$, then the subgroup $<H_1, K>$ is quasiconvex in $G$.*

Of course, Theorems 1 and 2 might be vacuously true, because there might be no subgroups $H_1$ and $K_1$ with the required properties. However, these theorems have interesting applications in the following frequently encountered case. Recall that a group $K$ is residually finite if for any finite set of non-trivial elements in $K$, there exists a finite index subgroup of $K$ which does not contain this set. The class of residually finite groups is very rich, it contains all finitely generated linear groups, and all fundamental groups of geometric 3-manifolds.

If $K$ is a quasiconvex subgroup of a negatively curved group $G$, then $K$ is finitely generated, so for any positive integer $n$, $K$ has only finitely many elements shorter than $n$. Hence if such $K$ is infinite and residually finite, it has an infinite family of distinct finite index subgroups $K_n$, such that $K_n$ does not contain non-trivial elements shorter than $n$. As $K_n$ is a finite index subgroup of a quasiconvex subgroup $K$, it is quasiconvex in $G$, hence we have the following results.

**Corollary 3.** *Let $H$ and $K$ be infinite quasiconvex subgroups of a negatively curved group $G$. Assume that $H$ and $K$ are residually finite. If the intersection of $H$ and $K$ is trivial, then there exist infinite families of distinct finite index subgroups $H_m$* of $H$ and $K_n$ of $K$, such that the subgroups $<H_m, K_n>$ are quasiconvex in $G$ and $<H_m, K_n> = H_m * K_n$.

**Corollary 4.** *Let $H$ and $K$ be infinite quasiconvex subgroups of a negatively curved group $G$. Assume that $H$ is malnormal in $G$ and residually finite. If the intersection of $H$ and $K$ is trivial, then there exists an infinite family of distinct finite index subgroups $H_m$ of $H$ with such that the subgroup $<H_m, K>$ is quasiconvex in $G$ and $<H_m, K> = H_m * K$.*

## Proofs of the Results.

Let $X$ be a set, let $X^* = \{x, x^{-1} | x \in X\}$, where $(x^{-1})^{-1} = x$ for $x \in X$. Denote the set of all words in $X^*$ by $W(X^*)$, and denote the equality of two words by " $\equiv$ ". Let $G$ be a group generated by the set $X^*$, and let $Cayley(G)$ be the Cayley graph of $G$ with respect to the generating set $X^*$. The set of vertices of $Cayley(G)$ is $G$, the set of edges of $Cayley(G)$ is $G \times X^*$, and the edge $(g, x)$ joins the vertex $g$ to $gx$.

**Definition.** The label of the path $p = (g, x_1)(gx_1, x_2)\cdots(gx_1x_2\cdots x_{n-1}, x_n)$
in $Cayley(G)$ is the word $Lab(p) \equiv x_1 \cdots x_n \in W(X^*)$. As usual, we identify the word $Lab(p)$ with the corresponding element in $G$. We denote the initial and the terminal vertices of $p$ by $\iota(p)$ and by $\tau(p)$ respectively, and the inverse of $p$ by $\bar{p}$. Denote the length of the path $p$ by $|p|$, where $|(g, x_1)(gx_1, x_2)\cdots(gx_1x_2\cdots x_{n-1}, x_n)| = n$.

A geodesic in the Cayley graph is a shortest path joining two vertices. A group $G$ is $\delta$-negatively curved if any side of any geodesic triangle in $Cayley(G)$ belongs to the $\delta$-neighborhood of the union of two other sides, (see [Gr] and [C-D-P]). Let $\lambda \leq 1, L > 0$ and $\epsilon > 0$. A path $p$ in the Cayley graph is a $(\lambda, \epsilon)$-quasigeodesic if for any subpath $p'$ of $p$ and for any geodesic $\gamma$ with the same endpoints as $p', |\gamma| > \lambda |p'| - \epsilon$.

A path $p$ is a local $(\lambda, \epsilon, L)$-quasigeodesic if for any subpath $p'$ of $p$ which is shorter than $L$ and for any geodesic $\gamma$ with the same endpoints as $p', |\gamma| > \lambda|p'| - \epsilon$, (cf. [C-D-P], p. 24).

Theorem 1.4 (p.25) of [C-D-P] (see also [Gr], p. 187) states that for any $\lambda_0 \leq 1$ and for any $\epsilon_0 > 0$ there exist constants $(L, \lambda, \epsilon)$ which depend only on $(\lambda_0, \epsilon_0)$ and $\delta$, such that any local $(\lambda_0, \epsilon_0, L)$-quasigeodesic in $G$ is a global $(\lambda, \epsilon)$-quasigeodesic in $G$.

Recall that $H$ is a $\mu$-quasiconvex subgroup of $G$ if any geodesic in $Cayley(G)$ which has its endpoints in $H$ belongs to the $\mu$-neighborhood of $H$.

*Proof of Theorem 1.* Let $H_1$ and $K_1$ be subgroups of $H$ and $K$ respectively, such that $H_1 \cap K_1 = H \cap K = G_0$. Consider an element $l$ of $<H_1, K_1>$ such that $l \notin G_0$. Then there exists a representation $l = h_1k_1 \cdots k_{m-1}h_m$, where $h_i \in H_1, k_i \in K_1, k_i$ and $h_i$ do not belong to $G_0, k_i$ and $h_i$ are geodesics in $G, h_1$ is a shortest representative of the coset $h_1G_0$, $h_m$ is a shortest representative of the coset $G_0h_m$, and for $1 < i < m, h_i$ is a shortest representative of the double coset $G_0h_iG_0$. (The elements $h_1$ or $h_m$ might be trivial). Let $p$ be the path in $Cayley(G)$ beginning at 1 with the decomposition of the following form: $p = p_1 q_1 \cdots q_{m-1} p_m$, where $Lab(p_i) \equiv h_i$ and $Lab(q_i) \equiv k_i$.

Let $A$ be the number of words in $G$ which are shorter than $2\mu + \delta$. As mentioned above, there exist constants $(L, \lambda, \epsilon)$ which depend only on $(\mu, A, \delta)$, such that any local $(1/3, (4\mu \cdot A + \delta), L)$-quasigeodesic in $G$ is a global $(\lambda, \epsilon)$-quasigeodesic in $G$.

Let $C_0 = max(L, \frac{\epsilon}{\lambda})$. We claim that if any elements in $H_1$ and in $K_1$ which are shorter than $C_0$ belong to $G_0$, then any path $p$, as above, is a $(\lambda, \epsilon)$-quasigeodesic in $G$. Indeed, it is enough to show that $p$ is a local $(1/3, (4\mu \cdot A + \delta), L)$-quasigeodesic in $G$. As $h_i \notin G_0$ and $k_i \notin G_0$, it follows that $|q_i| > C_0$ and $|p_i| > C_0$. As $L \leq C_0$, any subpath $t$ of $p$ with $|t| < L$ has a (unique) decomposition $t_1 t_2$, where without loss of generality, $t_1$ is a subpath of some $p_i$ and $t_2$ is a subpath of $q_i$. Let $t_3$ be a geodesic in $G$ joining the endpoints of $t$. We will show that $|t_3| \geq \frac{|t|}{3} - (4\mu \cdot A + \delta)$. If $|t_1| \leq \frac{|t|}{3}$, then $|t_3| \geq |t_2| - |t_1| \geq \frac{|t|}{3}$. So assume that $|t_1| > \frac{|t|}{3}$. As $t_1 t_2 t_3$ is a geodesic triangle in $Cayley(G)$, it is $\delta$-thin. Let $t_1'$ be the maximal subpath of $t_1$ which belongs to the $\delta$-neighborhood of $t_2$. It follows from Lemma 5 (below) that $|t_1'| \leq A \cdot 4\mu$, hence $|t_3| \geq |t_1| - |t_1'| - \delta \geq 1/3|t| - (A \cdot 4\mu + \delta)$.

So, indeed, $p$ is a local $(\frac{1}{3}, (4\mu \cdot A + \delta), L)$-quasigeodesic in $G$, hence it is a global $(\lambda, \epsilon)$-quasigeodesic in $G$. The definition of $C_0$ implies that if a $(\lambda, \epsilon)$-quasigeodesic $p$ in $G$ is longer than $C_0$, then $Lab(p) \neq 1$. As any element $l \in <H_1, K_1>$ which is not in $G_0$ has a representative $Lab(p)$ in $G$, as above, with $|p| > C_0$, it follows that $l \neq 1$, hence $<H_1, K_1> = H_1 *_{G_0} K_1$.

To prove the second part of the Theorem, for any element $l \in <H_1, K_1>$ which is not in $G_0$, consider a geodesic $\gamma$ in $G$ joining the endpoints of the path $p$, as above. As $G$ is negatively curved, there exists a constant $\alpha$ which depends only on $(\lambda, \epsilon)$ and $\delta$, such that $\gamma$ belongs to the $\alpha$-neighborhood of $p$. Assume that $H_1$ and $K_1$ are $\mu_1$-quasiconvex in $G$, then any vertex $v_i$ on $p_i$ is in the $\mu_1$-neighborhood of $h_1 k_1 \cdots k_{i-1} H_1$, and any vertex $w_i$ on $q_i$ is in the $\mu_1$-neighborhood of $h_1 k_1 \cdots h_i K_1$, so the path $p$ belongs to the $\mu_1$-neighborhood of the subgroup $<H_1, K_1>$. Hence $\gamma$ belongs to the $(\alpha + \mu_1)$-neighborhood of $<H_1, K_1>$.

If $l \in G_0$, then $l \in H_1$, so any geodesic $\gamma$ labeled with $l$ with $\iota(\gamma) = 1$, belongs to the $\mu_1$-neighborhood of $H_1$. It follows that the subgroup $<H_1, K_1>$ is $(\alpha + \mu_1)$-quasiconvex in $G$.

**Lemma 5.** *Using the notation of the proof of Theorem 1, $|t_1'| \leq A \cdot 4\mu$.*

*Proof.* To simplify notation, we drop the subscript $i$, so $t_1$ is a subpath of $p$, $t_2$ is a subpath of $q$, $Lab(p) = h$ and $Lab(q) = k$. As $H_1 < H$ and $K_1 < K$, we consider $h$ as an element of $H$ and $k$ as an element of $K$. Without loss of generality, assume that $q$ begins at 1, (so it ends at $k$), then $p$ begins at $h^{-1}$ and ends at 1. As $K$ and $H$ are $\mu$-quasiconvex in $G$, any vertex $v_i$ on $p$ is in the $\mu$-neighborhood of $H$, and any vertex $w_i$ on $q$ is in the $\mu$-neighborhood of $K$. Hence we can find vertices $v_1$ and $v_2$ in $t_1'$, $w_1$ and $w_2$ in $t_2$, $h'$ and $h''$ in $H$, and $k'$ and $k''$ in $K$ such that $|v_i, w_i| < \delta, |v_1, (h')^{(-1)}| < \mu, |v_2, (h'')^{(-1)}| < \mu, |w_1, k'| < \mu$ and $|w_2, k''| < \mu$. Then $|h'k'| < 2\mu + \delta$ and $|h''k''| < 2\mu + \delta$.

Assume that $|t_1'| > A \cdot 4\mu$. Then we can find vertices, as above which, in addition,

satisfy: $|v_2, v_1| > 4\mu$ and $h'k' = h''k''$. But then $(h'')^{(-1)}h' = k''(k')^{(-1)}$, so both products are in $G_0$. As $h$ is a shortest element in the double coset $G_0 h G_0$, it follows that $|h| \leq |h(h'')^{(-1)}h'|$. Let $r$ be a geodesic joining $(h'')^{(-1)}$ to $v_2$, let $s'$ be a subpath of $p$ joining $v_2$ to $1$ and let $s''$ be a subpath of $p$ joining $h^{-1}$ to $v_2$. Then $|h| = |p| = |s'| + |s''|$, and $|h(h'')^{(-1)}h'| \leq |s''| + |r| + |h'|$, hence $|s'| + |s''| \leq |s''| + |r| + |h'|$, so $|s'| + |r| \leq 2|r| + |h'|$. As $|h''| \leq |s'| + |r|$, and as $|r| \leq \mu$, it follows that $|h''| \leq 2\mu + |h'|$.

As $|v_2, v_1| > 4\mu$, the triangle inequality implies that $|h''| = |(h'')^{-1}| \geq |s'| - |r| = |1, v_1| + |v_1, v_2| - |r| \geq |1, v_1| + 4\mu - \mu = |1, v_1| + \mu + 2\mu$.

Now let $a$ be a geodesic joining $(h')^{-1}$ to $v_1$. As $|a| < \mu$, the triangle inequality implies that $|h'| = |(h')^{-1}| \leq |1, v_1| + |a| < |1, v_1| + \mu$. Hence, $|h''| > |h'| + 2\mu$, a contradiction. Therefore $|t'_1| \leq A \cdot 4\mu$.

We use the following property of malnormal quasiconvex subgroups of negatively curved groups, (cf. [Gi]). The original proof of this fact for the special case when $G$ is a free group is due to E. Rips ([G-R]).

Let $H$ be a subgroup of $G$, and let $G/H$ denote the set of right cosets of $H$ in $G$. The relative Cayley graph of $G$ with respect to $H$ is an oriented graph whose vertices are the cosets $G/H$, the set of edges is $(G/H) \times X^*$, such that an edge $(Hg, x)$ begins at the vertex $Hg$ and ends at the vertex $Hgx$. We denote it $Cayley(G, H)$. Note that for any path $p$ in $Cayley(G, H)$ if $\iota(p) = H \cdot 1$, then $\tau(p) = H \cdot Lab(p)$, so a path $p$ beginning at $H \cdot 1$ is closed, if and only if $Lab(p) \in H$.

**Lemma 6.** *Let $H$ be a malnormal $\mu$-quasiconvex subgroup of a finitely generated group $G$, and let $\delta$ be a non-negative constant. (The group $G$ does not have to be negatively curved). Let $\gamma_2 t \gamma_1$ be a path in $Cayley(G)$ such that $\gamma_1$ and $\gamma_2$ are geodesics in $Cayley(G)$, $Lab(\gamma_1) \in H, Lab(\gamma_2) \in H$ and $Lab(t) \notin H$. Let $m$ be the number of vertices in the ball of radius $\mu + 2\delta$ around $H \cdot 1$ in $Cayley(G, H)$, and let $M = m^2 + 1$. Then any subpath $\alpha$ of $\gamma_1$ which belongs to the $2\delta$-neighborhood of $\gamma_2$ in $Cayley(G)$ is shorter than $M$.*

*Proof of Lemma 6.* Assume that $|\alpha| \geq M$. Without loss of generality, $\iota(\gamma_2) = 1$. Let $\gamma'_1$ be a geodesic in $Cayley(G)$ beginning at $1$ such that $Lab(\gamma'_1) \equiv Lab(\gamma_1)$, and let $\alpha'$ be the subpath of $\gamma'_1$ such that $Lab(\alpha') \equiv Lab(\alpha)$. Let $\pi : Cayley(G) \to Cayley(G, H)$ be the projection map: $\pi(g) = Hg$ and $\pi(g, x) = (Hg, x)$.

Note that $\gamma_2$ and $\gamma'_1$ are geodesics in $Cayley(G)$ beginning at $1$, $Lab(\gamma'_1) \in H$ and $Lab(\gamma_2) \in H$, so as $H$ is $\mu$-quasiconvex in $G$, the projection $\pi$ maps $\gamma_2$ and $\gamma'_1$ into a ball of radius $\mu$ around $H \cdot 1$ in $Cayley(G, H)$. As $\alpha \subset N_{2\delta}(\gamma_2) \subset Cayley(G)$, it follows that $\pi(\alpha) \subset N_{2\delta+\mu}(H \cdot 1) \subset Cayley(G, H)$. Denote the vertices of $\pi(\alpha)$ by $v_1, \cdots, v_n$ and the vertices of $\pi(\alpha')$ by $w_1, \cdots, w_n$. As $n \geq M$, and $\{v_1, \cdots, v_n, w_1, \cdots w_n\} \subset N_{2\delta+\mu}(H \cdot 1)$, there exist $1 \leq i < j \leq n$ such that $(v_i, w_i) = (v_j, w_j)$. Let $p_1$ and $p_2$ be subpaths of $\pi(\gamma_1)$ connecting $\iota(\pi(\gamma_1))$ and $v_i$, and $v_i$ and $v_j$ respectively. Let $q_1$ and $q_2$ be subpaths of $\pi(\gamma'_1)$ connecting $H \cdot 1 = \iota(\pi(\gamma'_1))$ and $w_i$, and $w_i$ and $w_j$ respectively. Let $s = \pi(t)$ and let $\bar{p}$ denote the path $p$ with the opposite orientation. Then $sp_1 p_2 \bar{p}_1 \bar{s}$ and $q_1 q_2 \bar{q}_1$ are closed paths in $Cayley(G, H)$ beginning at $H \cdot 1$, so $Lab(sp_1 p_2 \bar{p}_1 \bar{s}) \in H$ and $Lab(q_1 q_2 \bar{q}_1) \in H$.

But $Lab(p_1) \equiv Lab(q_1)$ and $Lab(p_2) \equiv Lab(q_2)$, so $Lab(s)Lab(p_1p_2\bar{p}_1)Lab^{-1}(s) \in H$, and $Lab(p_1p_2\bar{p}_1) \in H$. Therefore the malnormality of $H$ in $G$ implies that $Lab(s) \in H$, contradicting the assumption that $Lab(s) = Lab(t) \notin H$. Hence $|\alpha| < M$.

*Proof of Theorem 2.* Let $H_1$ be a subgroup of $H$ such that $H_1 \cap K = H \cap K = G_0$. Let $l$ be an element of $<H_1, K>$ such that $l \notin G_0$. Then there exists a representation $l = h_1k_1 \cdots k_{m-1}h_m$, where $k_i$ and $h_i$ are as in the proof of Theorem 1. Let $p$ be a path in $Cayley(G)$ with a decomposition of the following form: $p = p_1q_1 \cdots q_{m-1}p_m$, where $Lab(p_i) \equiv h_i$ and $Lab(q_i) \equiv k_i$. Let $M$ be as in Lemma 6, and let $A$ be as in the proof of Theorem 1. As was mentioned above, there exist constants $(L', \lambda', \epsilon')$ which depend only on $(\mu, A, \delta)$, such that any local $(1/6, (4\mu \cdot A + \delta + M), L')$-quasigeodesic in $G$ is a global $(\lambda', \epsilon')$-quasigeodesic in $G$.

Let $C_1 = max(L', \frac{\epsilon'}{\lambda'})$. We claim that if all elements in $H_1$ which are shorter than $C_1$ belong to $G_0$, then any path $p$, as above, is a $(\lambda', \epsilon')$-quasigeodesic in $G$. Indeed, it is enough to show that $p$ is a local $(1/6, (4\mu \cdot A + \delta + M), L')$-quasigeodesic in $G$.

As $h_i \notin G_0$, it follows that $|p_i| > C_1$. As $L' < C_1$, any subpath $t$ of $p$ with $|t| < L'$ has a (unique) decomposition $t_1t_2t_3$, where $t_1$ and $t_3$ are subpaths of some $p_i$ and $p_{i+1}$, and $t_2$ is a subpath of $q_i$ (some of $t_i$ might be empty). Let $t_4$ be a geodesic in $G$ connecting the endpoints of $t$.

If $|t_2| > \frac{2|t|}{3}$, then $|t_1| + |t_3| \leq \frac{|t|}{3}$, so $|t_4| \geq |t_2| - (|t_1| + |t_3|) \geq \frac{2|t|}{3} - \frac{|t|}{3} = \frac{|t|}{3}$.

So assume that $|t_2| \leq \frac{2|t|}{3}$. Without loss of generality, assume that $|t_1| \geq |t_3|$, hence $|t_1| > \frac{|t|}{6}$. As $t_1t_2t_3t_4$ is a geodesic 4-gon in a $\delta$-negatively curved group $G$, there exists a decomposition $t_1 = s_2s_3s_4$ such that $s_2$ belongs to the $\delta$-neighborhood of $t_2$, $s_3$ belongs to the $2\delta$-neighborhood of $t_3$ and $s_4$ belongs to the $\delta$-neighborhood of $t_4$. According to Lemma 6, $|s_3| < M$ and according to Lemma 5, $|s_2| \leq 4\mu \cdot A$. But then $|t_4| + \delta \geq |s_4| = |t_1| - |s_2| - |s_3| \geq |t_1| - 4\mu \cdot A - M \geq \frac{|t|}{6} - 4\mu \cdot A - M$.

Hence $|t_4| \geq \frac{|t|}{6} - (M + \delta + 4\mu \cdot A)$, so the path $p$ is a local $(\frac{1}{6}, (M + \delta + 4\mu \cdot A), L')$-quasigeodesic in $G$, hence it is a $(\lambda', \epsilon')$-quasigeodesic in $G$. Then we conclude the argument as in the proof of Theorem 1.

## Acknowledgement

The problem solved in this paper appeared in M. Bestvina's problem list in geometric group theory. The author would like to thank Professors W. Neumann and J. Stallings for helpful conversations.

DEPARTMENT OF MATHEMATICS, CALTECH, PASADENA, CA 91125
*E-mail address*: `ritagtk @ cco.caltech.edu`